\documentclass[a4paper,11pt]{article}
\input epsf
\title{Matrices autosimilaires}
\author{Roland Bacher
}
\begin{document}
\maketitle

{\it R\'esum\'e: Cette note introduit une classe de matrices dont les 
d\'eterminants sont faciles \`a calculer. L'exemple le plus frappant
est obtenu en consid\'erant la matrice enti\`ere sym\'etrique $M(n)$
avec coefficients $M_{i,j}\in \{0,1\},\ 0\leq i,j<n$ 
d\'efinis en consid\'erant les coefficients  binomiaux modulo $2$ et en posant $M_{i,j}\equiv
{i+j \choose i}\pmod 2$. Le d\'eterminant $\hbox{det}(M(n))$ est alors 
\'etroitement reli\'e \`a la suite de Thue-Morse comptant les 
coefficients non-nuls d'entiers binaires.}

\section{D\'efinitions et r\'esultats}
Dans ce papier, la lettre $b$ d\'enotera toujours un entier $\geq 2$ et la
notation $n=\sum\nu_ib^i,\ 0\leq \nu_i\leq b-1$ d\'esignera un entier non-n\'egatif
$n$ \'ecrit en base $b$.

{\bf D\'efinition 1.1.} Une matrice (finie ou infinie)
$$M(n)=\left(M_{i,j}\right)_{0\leq i,j<n},\ n\in{\bf N}\cup\{\infty\}$$ 
(\`a coefficients $M_{i,j}$ dans le corps ou anneau unif\`ere commutatif favori du lecteur) 
est {\it $b-$autosimilaire} si on a $M_{0,0}=1$ et
$$M_{s,t}=\prod_i M_{\sigma_i,\tau_i}$$
pour tous $0\leq s=\sum \sigma_ib^i,\ t=\sum \tau_ib^i<n$.

Une matrice $M(n)$ qui est $b-$autosimilaire pour $b\leq n$ est enti\`erement d\'etermin\'ee
par la sous-matrice $\tilde M=M(b)$ de coefficients $\tilde M_{i,j}=M_{i,j},\ 0\leq i,j<b$.
On appelera $\tilde M$ la {\it matrice de d\'efinition} de  $M(n)$.

{\bf Exemple 1.2.} Pour $b=p$ un nombre premier, la matrice sym\'etrique $S$ et
la matrice triangulaire inf\'erieure $T$
$$S=\left(S_{i,j}={i+j\choose i}\pmod p\right)_{0\leq i,j}\quad \hbox{et}\quad
  T=\left(T_{i,j}={i\choose j}\pmod p\right)_{0\leq i,j}$$
\`a coefficients dans le corps fini ${\bf Z}/p{\bf Z}$ sont $p-$autosimilaires 
(voir section 2).

{\bf Remarque 1.3.}
Une autre raison de s'int\'eresser aux matrices $b-$autosimilaires
provient du produit tensoriel:
Soit $(\tilde B_{i,j})_{0\leq i,j<b}$ une matrice \`a coefficients dans un corps $K$
d\'efinissant un endomorphisme de $V=K^b$.
La matrice $B(b^d)$ (qui n'est $b-$autosimilaire que pour $\tilde B_{0,0}=1$)
d\'efinie par
$$B_{s,t}=\prod_{i=0}^{d-1} 
\tilde B_{\sigma_i,\tau_i},\ 0\leq s=\sum_{i=0}^{d-1}
\sigma_i b^i,\ t=\sum_{i=0}^{d-1}\tau_i b^i<b^d$$
d\'efinit alors un endomorphisme diagonal 
$$\begin{array}{cccc}
\displaystyle \tilde B^{\otimes^d}:&\displaystyle V^{\otimes^d}&\displaystyle \longrightarrow
&\displaystyle V^{\otimes^d}\cr
&\displaystyle v_0\otimes\cdots \otimes v_{i}\otimes\cdots
\otimes v_{d-1}&\displaystyle \longmapsto
&\displaystyle \tilde Bv_0\otimes\cdots \otimes \tilde B v_i\otimes \cdots
\otimes \tilde Bv_{d-1}\end{array}$$
o\`u $v_0,\dots,v_{d-1}\in V=K^b$. La normalisation $B_{0,0}=1$ permet
alors de s'affranchir des puissances enti\`eres $b^d$ de $b$ et d\'efinit
un analogue d'une telle matrice pour toute dimension enti\`ere ainsi que pour
une dimension $\infty$ d\'enombrable.

Une matrice $b-$autosimilaire $M(n),\ n\in\{b,b+1,\dots\}\cup\{\infty\}$, est 
{\it non-d\'eg\'en\'er\'ee}
si les $b$ matrices $M(k)=(M_{i,j})_{0\leq i,j<k},\ k=1,\dots,b$ sont toutes inversibles.
On pose alors $d(0)=1$ et 
$$d(k)=\hbox{det}(M(k+1))/\hbox{det}(M(k)),\ k=1\dots,b-1\ .$$

En particulier, une matrice $b-$autosimilaire triangulaire $T(\infty)$
est non-d\'eg\'en\'er\'ee si et seulement
si sa matrice de d\'efinition $\tilde T$ est
inversible. Les nombres $d(0),\dots,d(b-1)$ introduits ci-dessus ne
sont alors rien d'autre que les coefficients diagonaux de $\tilde T$.

{\bf Th\'eor\`eme 1.4.} {\sl Soit $b\geq 2$ un entier et $M=M(n)$ une
matrice $b-$autosimilaire non-d\'eg\'en\'er\'ee.

\ \ (i) On a une factorisation unique
$$M=LDU$$
o\`u $L$ est une matrice $b-$autosimilaire,
unipotente ($L_{i,i}=1,\ 0\leq i<n$) triangulaire inf\'erieure, $D$ est 
$b-$autosimilaire diagonale
et $U$ est $b-$autosimilaire, unipotente triangulaire sup\'erieure.

\ \ (ii) On a $\tilde D_{i,i}=d(i),\ 0\leq i<b$ et
$$\hbox{det}(M(n))=\hbox{det}(D(n))=\prod_{k=0,\ k=
\sum \kappa_ib^i}^{n-1}\qquad\prod_{i} d(\kappa_i)\not=0$$
pour tout $n\in{\bf N}$.

\ \ (iii) L'ensemble des matrices (infinies) triangulaires
 inf\'erieures (ou sup\'erieures)
$b-$autosimilaires et non-d\'eg\'en\'er\'ees est un
groupe: Si $R$ et $T$ sont deux telles matrices d\'efinies par 
$\tilde R$ 
et $\tilde T$, alors $RT$ est d\'efinie par $\tilde R\tilde T$.}

{\bf Remarques 1.5.} (i) Une matrice $M(n),\ n\in{\bf N}$ qui est $b-$autosimilaire
et non-d\'eg\'en\'er\'ee s'inverse donc relativement facilement:
$$\left(M(n)\right)^{-1}=\left(U(n)\right)^{-1}\ \left(D(n)\right)^{-1}\ 
\left(L(n)\right)^{-1}$$
avec $U,D,L$ comme dans l'assertion (i) trois matrices inversibles triangulaires.
Leurs inverses sont donc $b-$autosimilaires et peuvent se calculer \`a 
partir des inverses $\tilde U^{-1},\tilde D^{-1},\tilde L^{-1}$ de leurs matrices
de d\'efinition.

\ \ (ii)  Le th\'eor\`eme 1.4 permet tr\`es souvent de calculer 
$\hbox{det}(M(n))$ m\^eme si la matrice $b-$autosimilaire $M=M(n)$
(avec $n\in{\bf N}$) n'est pas non-d\'eg\'en\'er\'ee: Il suffit de 
consid\'erer la matrice $b-$autosimilaire associ\'ee \`a une perturbation
g\'en\'erique $\tilde D+t\tilde A$ (avec $\tilde A_{0,0}=0)$) et d'\'evaluer
le r\'esultat en $t=0$. L'exemple 2.4 de la section suivante illustrera 
ce proc\'ed\'e. 

\ \ (iii) Si la matrice $\tilde B$ consid\'er\'ee dans la remarque 1.3 est une 
matrice inversible, alors la preuve du th\'eor\`eme montre 
que la matrice $B(b^d)$ est inversible et que son inverse est $C(b^d)$
o\`u $\tilde B\tilde C=\hbox{id}$. 

{\bf Preuve du th\'eor\`eme 1.4.} 
Soit $\tilde M$ la matrice de d\'efinition de la matrice $b-$autosimilaire
$M(n)$. Comme $M(n)$ est non-d\'eg\'en\'er\'ee, on a une factorisation unique
(donn\'ee par le proc\'ed\'e d'orthogonalisation de Gram-Schmitt)
$$\tilde M=\tilde L\tilde D\tilde U$$
avec $\tilde L$ unipotent triangulaire inf\'erieure, $\tilde D$ diagonale et 
$\tilde U$ unipotent
triangulaire sup\'erieure. Comme $\tilde L_{0,0}=\tilde D_{0,0}=\tilde U_{0,0}=1$, les matrices
$\tilde L,\tilde D$ et $\tilde U$ d\'efinissent des 
matrices $b-$autosimilaires $L(l),D(l)$ et $U(l)$ pour tout $l\in{\bf N}\cup\{\infty\}$ 
et un petit calcul montre que ces matrices sont respectivement unipotente
triangulaire inf\'erieure, diagonale et unipotente triangulaire sup\'erieure.

On a alors (avec $L=L(\infty),D=D(\infty),U=(\infty)$ et $M=M(\infty)$) pour
$0\leq s=\sum\sigma_ib^i,\ t=\sum \tau_ib^i$
$$\begin{array}{ll}
\displaystyle \left(LDU\right)_{s,t}&
\displaystyle =\sum_k L_{s,k}D_{k,k}U_{k,t}\cr
&\displaystyle =\sum_{k=\sum \kappa_i b^i}\quad  \prod_i \tilde L_{\sigma_i,\kappa_i}
\tilde D_{\kappa_i,\kappa_i}\tilde U_{\kappa_i,\tau_i}\cr
&\displaystyle =\prod_i \sum_{\kappa_i=0}^{b-1} \tilde L_{\sigma_i,\kappa_i}
\tilde D_{\kappa_i,\kappa_i}\tilde U_{\kappa_i,\tau_i}\cr
&\displaystyle =\prod_i \tilde M_{\sigma_i,\tau_i}
=M_{s,t}\end{array}$$
ce qui montre l'assertion (i) (modulo l'unicit\'e, laiss\'ee au lecteur).

Les \'egalit\'es
$$\hbox{det}(M(l))=\hbox{det}(D(l))\quad \hbox{pour }l=1,\dots,b$$
impliquent par r\'ecurrence sur $0\leq i<b$ qu'on a $\tilde D_{i,i}=D_{i,i}=d(i)$.
On a donc
$$\hbox{det}(M(n))=\hbox{det}(D(n)=\prod_{k=0}^{n-1} D_{k,k}$$
$$=\prod_{k=0,\ k=\sum\kappa_ib^i}^{n-1}\qquad \prod_{i}
 \tilde D_{\kappa_i,\kappa_i}$$
ce qui prouve l'assertion (ii).

L'assertion (iii) r\'esulte d'un calcul similaire \`a celui utilis\'e
dans la preuve de l'assertion (i), \`a savoir:
$$\begin{array}{ll}
\displaystyle \left(RT\right)_{s,t}&
\displaystyle =\sum_k R_{s,k}T_{k,t}\cr
&\displaystyle =\sum_{k=\sum \kappa_i b^i}\quad  \prod_i \tilde R_{\sigma_i,\kappa_i}
\tilde T_{\kappa_i,\tau_i}\cr
&\displaystyle =\prod_i \sum_{\kappa_i=0}^{\sigma_i} \tilde R_{\sigma_i,\kappa_i}
\tilde T_{\kappa_i,\tau_i}\cr
&\displaystyle =\prod_i \left(\tilde R\tilde
  T\right)_{\sigma_i,\tau_i}
\end{array}$$
pour $R$ et $T$ deux matrices triangulaires inf\'erieures 
$b-$autosimilaires. Il faut \'egalement v\'erifier la condition de 
normalisation $(RT)_{0,0}=1$ qui r\'esulte \'evidemment du fait
que $R$ et $T$ sont tous les deux triangulaires inf\'erieures.\hfill QED

\section{Exemples li\'es au triangle de Pascal}

Fixons un nombre premier $p$ et consid\'erons la r\'eduction $\pmod p$ du
triangle de Pascal form\'e des coefficients binomiaux ${n\choose i},\ n,i\in{\bf N}$ 
d\'efinis par
$$(1+x)^n=\sum_i {n\choose i}x^i\in {\bf Z}[x]\ .$$

L'existence de l'automorphisme de Frobenius $\varphi:\alpha\longmapsto \alpha^p$
sur un corps de caract\'eristique $p$ montre qu'on a pour 
$n=\sum\nu_ip^i$ 
$$(1+x)^n\equiv \prod_i \left(1+x^{p^i}\right)^{\nu_i}\pmod p$$
ce qui implique pour $k=\sum\kappa_ip^i$ la factorisation
$${n\choose k}=\prod {\nu_i\choose \kappa_i}\pmod p\leqno{\hbox{(Eq. 1)}}$$
et prouve que les matrices de l'exemple 1.2 sont bien $p-$autosimilaires.

L'identit\'e triviale $(1+x)^{s+t}=(1+x)^s\ (1+x)^t$ implique l'\'egalit\'e
$${s+t\choose s}=\sum_k {s\choose s-k}{t\choose k}=\sum_k {s\choose k}{t\choose k}
$$
et fournit la factorisation $S=S(n)=L(n)L(n)^t$ de la matrice sym\'etrique enti\`ere
$S$ avec coefficients $S_{i,j}={i+j\choose i},\ 0\leq i,j<n$ comme produit
de la matrice unipotente triangulaire inf\'erieure $L$ avec coefficients $L_{i,j}={i\choose j},
\ 0\leq i,j<n$ et de sa transpos\'e. La r\'eduction $\pmod p$ de cette factorisation
illustre l'assertion (i) du Th\'eor\`eme 1.4 (avec $D=\hbox{id}$ et $U=L^t$)
pour les matrices de l'exemple 1.2.

Consid\'erons la matrice 
$$A=A(n)=\left(\begin{array} {cccccc}
0&0&0&\dots\cr
1&0&\cr
0&2&\cr
\vdots&&\ddots\end{array}\right)$$
d\'efinie par
$$A_{i,j}=\left\{\begin{array}{ll}
\displaystyle i\qquad&\displaystyle \hbox{si }i=j+1\ ,\cr
\displaystyle 0\qquad&\displaystyle \hbox{sinon}\end{array}\right.$$
pour $0\leq i,j<n$. Un calcul facile montre qu'on a
$$L(n)=\hbox{ exp }A(n)=\sum_{k=0}^\infty \frac{A^k}{k!}$$
(o\`u $L(n)$, donn\'e par $L_{i,j}={i\choose j},\ 0\leq i,j<n$, est
la matrice unipotente triangulaire inf\'erieure introduite ci-dessus)
ce qui implique que la matrice $R=L^{-1}=\hbox{exp } (-A)$ est donn\'e par
$$R_{i,j}=(-1)^{i+j} {i\choose j},\ 0\leq i,j$$
(ceci se d\'emontre \'evidemment aussi ais\'ement en calculant $LR$). Le
calcul de la matrice $S(n)^{-1}=R(n)^t \ R(n), n\in {\bf N}$ (dans l'anneau ${\bf Z}$ des entiers ou dans
le corps fini ${\bf Z}/p{\bf Z}$) est donc tout \`a fait explicite.

Soit $p$ un premier que nous fixons et soit 
$a\in{\bf Z}/p{\bf Z}$ un \'el\'ement du corps fini \`a
$p$ \'el\'ements. Rappelons que le symbole de
Legendre $\left(\frac{a}{p}\right)$ est la fonction d\'efinie par
$$\left(\frac{x}{p}\right)=\left\{
\begin{array}{cll}
\displaystyle 0&\qquad&\displaystyle \hbox{si } a=0\ ,\cr
\displaystyle 1&&\displaystyle \hbox{si } a\not=0\hbox{ est un carr\'e de }
{\bf Z}/p{\bf Z}\ ,\cr
\displaystyle -1&&\displaystyle \hbox{si } a\hbox{ n'est pas
un carr\'e dans }{\bf Z}/p{\bf Z}\ .\end{array}\right.$$

Posons $\psi(a)=\left(\frac{a}{p}\right)$ et consid\'erons la
matrice enti\`ere sym\'etrique $M(n)$ avec coefficients
$$M_{i,j}=\psi({i+j\choose i}\pmod p)\ ,0\leq i,j<n\ .$$
Les propri\'et\'es du symbole de Legendre et (Eq. 1) montrent que $M$ est 
$p-$autosimilaire.

{\bf Exemple 2.1.}
Pour $p=2$ la matrice $M(\infty)$ n'est rien d'autre que la r\'eduction
modulo $2$ du triangle de Pascal (sous forme de matrice sym\'etrique)
avec coefficients dans $\{0,1\}$. On a la factorisation
$$\tilde M=\left(\begin{array}{rrrrrrrrrrrr}
1&1\cr
1&0\end{array}\right)=\left(\begin{array}{rrrrrrrrrrrr}
1&0\cr
1&1\end{array}\right)
\left(\begin{array}{rrrrrrrrrrrr}
1&0\cr
0&-1\end{array}\right)
\left(\begin{array}{rrrrrrrrrrrr}
1&1\cr
0&1\end{array}\right)$$
de la matrice de d\'efinition $\tilde M$ et le Th\'eor\`eme 1.4
montre alors qu'on a 
$$\hbox{det}(M(n))=\prod_{\sum\kappa_i2^i=0}^{n-1}(-1)^{\sum \kappa_i}
=\prod_{k=0}^{n-1}(-1)^{s_k}$$
o\`u la suite $s_k$ (d\'efinie r\'ecursivement par $s_0=0$, 
$s_{2k}=s_{k}$, et $s_{2k+1}=1-s_k$) est la fameuse suite de Thue-Morse
comptant $\pmod 2$ le nombre de digits non-nuls d'un entier binaire
(cf. {\bf [AS]}). Plus pr\'ecis\'ement, on montre facilement qu'on a 
$$\hbox{det}(M(2n))=(-1)^n\quad\hbox{ et }\quad
\hbox{det}(M(2n+1))=(-1)^{n+s_n}\ .$$

Le fait que la matrice $\tilde L \tilde D=\left(\begin{array}{rr}
1&0\cr 1&-1\end{array}\right)$ (et donc $L(\infty)D(\infty)$) soit
d'ordre $2$ est \'equivalente \`a la d\'efinition r\'ecursive
suivante de la suite de Thue-Morse: $s_0=0$ et $s_n\in\{0,1\}$ tel que
$$\sum_k\left({n\choose k}\pmod 2\right)(-1)^{s_k}=0$$
pour tout $n\geq 1$ (ici, $\left({n\choose k}\pmod 2\right)\in\{0,1\}$
et l'\'egalit\'e est sur ${\bf Z}$). Ceci peut se g\'en\'eraliser
comme suit: Le vecteur $\mu=(1,-1,-1,1,\dots,(-1)^{s_n},\dots)$ associ\'e
\`a la suite de Thue-Morse v\'erifie
$$L\ \mu^t=\left(\begin{array}{c}1\cr 0\cr\vdots\end{array}\right)\ ,
\qquad \hbox{i.e.}\quad \mu^t=L^{-1}
\left(\begin{array}{c}1\cr 0\cr\vdots\end{array}\right)$$
o\`u $L$ est la matrice $2-$autosimilaire triangulaire inf\'erieure
d\'efinie par $\tilde
L=\left( \begin{array}{cc}1&0\cr
1&1\end{array}\right)$. L'assertion (iii) du Th\'eor\`eme 1.4 montre
donc
que $L^{-1}$ est $2-$autosimilaire, d\'efini par $\tilde L^{-1}=
\left( \begin{array}{rr}1&0\cr
-1&1\end{array}\right)$. En rempla\c cant $L$ ci-dessus par la
matrice $b-$autosimilaire triangulaire inversible 
d\'efinie par $\tilde L_{i,j}={i\choose j},\ 0\leq i,j<b$, d'inverse
$\tilde R$ avec coefficients $\tilde R_{i,j}=(-1)^{i+j}
{i\choose j},\ 0\leq i,j<b$, on a 
$$L\ \mu^t=\left(\begin{array}{c}1\cr 0\cr\vdots\end{array}\right)\ ,
\qquad \hbox{i.e.}\quad \mu^t=R
\left(\begin{array}{c}1\cr 0\cr\vdots\end{array}\right)$$
pour le vecteur
$\mu=(\mu_0,\mu_1,\dots)=(1,-1,\dots,\mu_n=(-1)^{\sum \nu_i},\dots)$
avec $n=\sum_i \nu_i b^i$. En particulier, si $b$ est impair, on a 
simplement $\mu_n=(-1)^n$ ce qui se traduit par l'\'egalit\'e
$$\sum_{k\equiv 0\pmod 2}\quad \prod_i {\nu_i\choose \kappa_i}= 
\sum_{k\equiv 1\pmod 2}\quad\prod_i {\nu_i\choose \kappa_i}=
2^{\left(\sum\nu_i\right)-1}$$
pour $n=\sum \nu_ib^i \geq 1$ et $k=\sum \kappa_i b^i$ avec $b\geq 3$
un entier impair. On peut \'evidemment \'egalement consid\'erer 
d'autres matrices $b-$autosimilaires triangulaires inf\'erieures
$\Lambda$ \`a 
la place de la matrice $L$ discut\'ee ci-dessus et calculer la suite $\mu$
associ\'ee (qui n'est rien d'autre que la suite des coordonn\'ees de la
premi\`ere colonne de la matrice $b-$autosimilaire $\Lambda^{-1}$).

Revenons maintenant au premier $p=2$ et \`a notre matrice sym\'etrique
$2-$autosimilaire $M$ d\'efinie par $\tilde M_{i,j}={i+j\choose
  i}\pmod 2\in\{0,1\}$.
On peut montrer assez facilement que
 l'inverse $M(n)^{-1}$ de la matrice $M(n)$ (pour $n\in {\bf N}$)
n'a que des coefficients dans $\{-1,0,1\}$.  

Signalons finalement que les d\'eterminants des
sous-matrices de $M(\infty)$ obtenu en prenant
$l$ lignes et colonnes cons\'ecutives (et ayant donc comme 
coefficients $M_{i+s,j+t},\ 0\leq i,j<l$ pour $s,t\in{\bf N}$ fix\'e)
ne semblent prendre que les valeurs $0$ ou $\pm 1$.

{\bf Exemple 2.2.}
Pour le premier $p=3$, la matrice $M(\infty)$ est la r\'eduction 
du triangle de Pascal $\pmod 3$ \`a valeurs dans $\{0,\pm 1\}$.
La factorisation
$$\tilde M=\left(\begin{array}{rrrrrrrrrrrr}
1&1&1\cr
1&-1&0\cr
1&0&0\end{array}\right)=\left(\begin{array}{rrrrrrrrrrrr}
1&0&0\cr
1&1&0\cr
1&\frac{1}{2}&1
\end{array}\right)
\left(\begin{array}{rrrrrrrrrrrr}
1&0&0\cr
&-2&0\cr
0&0&-\frac{1}{2}
\end{array}\right)
\left(\begin{array}{rrrrrrrrrrrr}
1&1&1\cr
0&1&\frac{1}{2}\cr
0&0&1
\end{array}\right)$$
montre que le d\'eterminant
$$\hbox{det}(M(n))=\prod_{\sum\kappa_i3^i=0}^{n-1}\prod_i d(\kappa_i)$$
(avec $d(0)=1,\ d(1)=(-2)$ et $d(2)=-\frac{1}{2}$) est une
puissance de $(-2)$.

{\bf Example 2.3.}
Pour $p=5$ on a $\tilde M=\tilde L\tilde D\left(\tilde L\right)^t$ avec
$$\tilde M=\left(\begin{array}{rrrrrrrrrrrr}
1&1&1&1&1\cr
1&-1&-1&1&0\cr
1&-1&1&0&0\cr
1&1&0&0&0\cr
1&0&0&0&0
\end{array}\right)\ ,\qquad \tilde L\left(\begin{array}{rrrrrrrrrrrr}
1&0&0&0&0\cr
1&1&0&0&0\cr
1&1&1&0&0\cr
1&0&-\frac{1}{2}&1&0\cr
1&\frac{1}{2}&0&\frac{2}{3}&1
\end{array}\right)$$
 et $\tilde D$ une matrice diagonale avec coefficients diagonaux
$$d(0)=1,\ d(1)=-2,\ d(2)=2,\ d(3)=-\frac{3}{2},\ d(4)=\frac{1}{6}\ .$$
Le d\'eterminant
$$\hbox{det}(M(n))=\prod_{\sum \kappa_i5^i=0}^{n-1}\quad
\prod_i d(\kappa_i)$$
est donc un entier de la forme $\pm 2^{\alpha}3^{\beta}$ avec
$\alpha=\alpha(n)$ et $\beta=\beta(n)$ des entiers naturels convenables.

{\bf Example 2.4.}
Pour le premier $p=7$ la matrice $M(\infty)$ n'est plus non-d\'eg\'en\'er\'ee.
On peut la perturber en consid\'erant la matrice $p-$autosimilaire
d\'efinie par
$$\tilde M=\left(\begin{array}{rrrrrrrrrrrr}
1&1&1&1&1&1&1\cr
1&t+1&-1&1&-1&-1&0\cr
1&-1&-1&-1&1&0&0\cr
1&1&-1&-1&0&0&0\cr
1&-1&1&0&0&0&0\cr
1&-1&0&0&0&0&0\cr
1&0&0&0&0&0&0
\end{array}\right)$$
qui est maintenant non-d\'eg\'er\'ee sur ${\bf Z}[t]$ et qui redonne la 
matrice initiale en posant $t=0$.

La matrice triangulaire $\tilde L$ est alors donn\'ee par
$$\tilde L=\left(\begin{array}{rrrrrrrrrrrr}
1&0&0&0&0&0&0\cr
1&1&0&0&0&0&0\cr
1&\frac{-2}{t}&1&0&0&0&0\cr
1&0&\frac{t}{t+2}&1&0&0&0\cr
1&\frac{-2}{t}&\frac{2}{t+2}&\frac{t-2}{4}&1&0&0\cr
1&\frac{-2}{t}&\frac{t+4}{2t+4}&\frac{-1}{2}&\frac{-6}{t-10}&1&0\cr
1&\frac{-1}{t}&\frac{1}{2}&0&\frac{-4}{t-10}&\frac{t+2}{t+8}&1\end{array}
\right)$$
et la matrice diagonale $\tilde D$ correspondante est donn\'ee par
$$d(0)=1,\ d(1)=t,\ d(2)=\frac{-2t-4}{t},\ d(3)=\frac{-4}{t+2},\ $$
$$d(4)=
\frac{t-10}{4},\ d(5)=\frac{-t-8}{2t-20},\ d(6)=\frac{-1}{t+8}\ .$$

Le d\'eterminant $\hbox{det}(M(n))$ se calcule maintenant en \'evaluant 
la fraction rationelle (qui est en fait toujours un \'el\'ement de
${\bf Z}[t]$)
$$\hbox{det}(M(n))=\prod_{\sum \kappa_i7^i=0}^{n-1}\quad
\prod_i d(\kappa_i)$$
en $t=0$. Il est soit $0$ soit de la forme $\pm 2^\alpha\ 5^\beta$.
On peut montrer qu'il est nul si et seulement si $m=\sum\mu_i 7^i=n-1$
contient le chiffre $\mu_i=1$ ou si $m$ contient le chiffre $\mu_i=2$ suivi
d'autre chose que de $666\dots$ (i.e. $\hbox{det}(M(1+\sum \mu_i7^i))=0
\Longleftrightarrow$ ou bien il existe $i$ avec $\mu_i=1$ ou bien il existe
$i>j$ avec $\mu_i=2$ et $\mu_j\not=6$).

Je remercie Jean-Paul Allouche pour ses commentaires et son int\'er\^et.

\vskip.5cm
$\hbox{\bf [AS]}$ J-P.Allouche, J.Shallit, {\it The ubiquitous Prouhet-Thue-Morse
sequence}, Proceedings of SETA 98 (C.Ding, T.Helleseth, H.Niederreiter, editors),
Springer (1999).

\vskip.8cm
Roland Bacher

Institut Fourier, UMR 5582

Laboratoire de Math\'ematiques, BP 74, 38402 St. Martin d'H\`eres Cedex,
France

{\it E-mail}: Roland.Bacher@ujf-grenoble.fr

\end{document}